\documentstyle{article}
\def\Bbb{\bf}
\def\bea{\begin{eqnarray*}}
\def\eea{\end{eqnarray*}}
\newtheorem{main}{Theorem}

\newtheorem{defn}{Definition}
\newtheorem{thm}{Theorem}
\newenvironment{proof}{\medskip {\bf Proof.}}{\hfill \rule{.5em}{1em}
\\}
\begin{document}
\sloppy
\title{4-Manifolds without Einstein Metrics}

\author{Claude LeBrun\thanks{Supported
in part by  NSF grant DMS-9505744.}
\\
SUNY Stony
 Brook
  }

\date{November, 1995}
\maketitle

\begin{abstract}
It is shown that there are  infinitely many  compact
orientable smooth 4-manifolds which do not
admit Einstein metrics,  but nevertheless satisfy
the strict Hitchin-Thorpe inequality
$2\chi > 3|\tau |$.
 The examples in question arise as
non-minimal complex algebraic surfaces of general type,
and the method of proof stems from Seiberg-Witten theory.
 \end{abstract}

\section{Introduction}

A smooth Riemannian metric $g$ is said to be {\em  Einstein}
if its Ricci curvature $r$ is a constant multiple of
the metric:
$$r=\lambda g .$$
Not every   4-manifold admits such metrics.
A necessary condition for the existence of an Einstein metric
on a compact oriented 4-manifold is that the
Hitchin-Thorpe inequality $2\chi  (M) \geq  3 |\tau (M)|$
must hold \cite{bes}. Moreover, equality can  hold only if $M$
manifold is finitely covered by a torus or K3 surface. We will
  say that  $M$  satisfies the {\em   strict
 Hitchin-Thorpe inequality} if $2\chi  (M) >  3 |\tau (M)|$.

The purpose of this note is to prove the following result:

\begin{main}
There  are  infinitely many  compact simply-connected
 smooth 4-manifolds which
do not
admit Einstein metrics,
but nevertheless
 satisfy the strict
Hitchin-Thorpe inequality.
\end{main}

The examples we shall consider arise as non-minimal
complex surfaces of general type. The proof hinges on
scalar curvature estimates that come from Seiberg-Witten
theory.

\section{Scalar Curvature and Topology}

In this section, we will develop a certain lower bound for the
$L^2$-norm of the scalar curvature of all Riemannian
metrics on a non-minimal complex surfaces of general type.
Let us begin by reviewing some definitions and results.

\begin{defn}
Let $M$ be a smooth compact oriented 4-manifold.
A {\em polarization} of $M$ is a linear subspace
$H^+\subset H^2(M, {\Bbb  R})$  on which
the restriction of the intersection form
 is positive-definite, and which is a maximal
 subspace with this property.
\end{defn}

The example of interest  is the following:
let $g$ be a Riemannian metric, and let
$H^+(g)$ be the space of harmonic self-dual
2-forms with respect to $g$. Then $H^+(g)$
is a polarization. If $H^+$ is a given polarization,
and if $H^+(g)=H^+$, we will say that $g$ is
adapted to $H^+$.

If $J$ is an orientation-compatible almost-complex structure on
$M$,  $J$ induces a spin$^c$-structure $c$ on $M$, and for
every metric $g$ one then has a pair of rank-2 complex
vector bundles $V_\pm$ which formally
satisfy
$$V_\pm = {\Bbb S}_\pm \otimes L^{1/2} , $$
where ${\Bbb S}_{\pm}$ are the left- and right-handed
spinor bundles of $g$, and  $L$ is the anti-canonical
line bundle of $J$. For each unitary connection $\theta$ on $L$,
we  have a Dirac operator
$D_{\theta}: C^{\infty}(V_+)\to C^{\infty}(V_-)$,
and one can then consider the Seiberg-Witten equations
\cite{witten}
\bea D_{\theta}\Phi &=&0 \\
 F_{\theta}^+&=&i \sigma(\Phi) \eea
for an unknown section $\Phi$ of $V_+$
and an unknown unitary connection $\theta$.
 Suppose that $H^+$ is a polarization such that
the   orthogonal projection $c^+$ of $c_1(L)$ into
$H^+$ is non-zero. Let $g$ be any $H^+$-adapted metric,
and consider the moduli space of solutions
of a generic perturbation of the Seiberg-Witten
equations modulo gauge equivalence. This moduli space consists
of a finite number of oriented points, and
the Seiberg-Witten invariant $n_c(M,H^+)$
is defined to be the number of points
in moduli space, counted with signs.
This is independent of all choices. Indeed,
if $b^+(M) > 1$, it is even independent of $H^+$.

A Weitzenb\"ock argument yields
the following curvature estimate \cite{leb2}:

\begin{thm}\label{est}
Let $(M,H^+,c)$ be a  smooth compact
oriented   polarized 4-manifold with spin$^c$ structure
such that  $n_c(M,H^+)\neq 0$. If $c_1(L)\in H^2(M, {\Bbb R})$ is
the anti-canonical class of this structure,
  let $c_1^+\neq 0$ be its orthogonal projection to
$H^+$ with respect to the intersection form.  Then  every
$H^+$-adapted Riemannian metric $g$ satisfies
$$ \int_M s^2~d\mu \geq 32\pi^2(c_1^+)^2 , $$
with equality iff $g$ is  K\"ahler with respect to a
$c$-compatible
complex structure and
 has constant negative
scalar curvature.
 \end{thm}

Now if $(M,J)$ is a complex surface of K\"ahler type with $b^+>1$, and
if $c$ is the
spin$^c$ structure induced by $J$, then $n_c(M,H^+)=n_c(M)=1$.
For complex surfaces with $b^+=1$, the picture is more complicated,
but can be summarized as follows.  The set of classes $\alpha \in
H^2(M,{\Bbb R})$
with $\alpha^2:=\alpha \cdot \alpha > 0$ consists of two connected
components.
One component contains the K\"ahler classes of  all K\"ahler metrics
on $M$;  let us call the elements of this component
{\em future pointing}, and  the elements of the other
{\em past-pointing}.
Then $n_c(M,H^+)=1$ if $c^+$ is past-pointing,
and that $n_c(M,H^+)=0$ if $c^+$ is future pointing \cite{FM,KM}.

We  now come to the   technical  heart of the article:

\begin{thm}\label{best}
Let $X$ be a minimal complex algebraic surface
of general type,
and let $M= X\#k\overline{\Bbb CP}_2$ be obtained from
$X$ by blowing up $k >0$ points.
Then any Riemannian metric on $M$ satisfies
$$\int_M s^2 d\mu > 32\pi^2 (2 {\chi}+3 {\tau} +k),$$
where $ {\chi}$ and ${\tau}$ are respectively
the Euler characteristic and signature of   $M$.
 \end{thm}

\begin{proof}
Let us  think of $M$ concretely as obtained from
$X$ by blowing up $k$ distinct points
$p_1$, \ldots , $p_k$, so that $M$ comes
equipped with an integrable complex structure $J$.
The key observation \cite{FM} is that instead of merely
considering this complex
structure alone, it is natural to
 consider $2^k$ distinct complex structures,
each of which is the pull-back  of $J$  via a diffeomorphism
$M\to M$. To this end, choose a  biholomorphism between
a neighborhood of    $p_j\in X$ and the unit ball in ${\Bbb C}^2={\Bbb
R}^4$.  Let
 $\psi_j: X\to X$ be the identity outside this neighborhood,  and act
by
$$\left[ \begin{array}{cccc}
1&0&0&0\\ 0&\cos \pi u (r)& 0& -\sin \pi u(r) \\
0&0&1&0 \\ 0 & \sin \pi u(r) & 0& \cos \pi u(r)
\end{array} \right]  $$
on the ball itself;
here $r$ is the distance from the origin
in ${\Bbb R}^4$,  and the smooth function $u$
satisfies $u(r) \equiv 1$ for $r\leq \frac{1}{3}$  and
 $u(r) \equiv 0$ for
$r\geq \frac{2}{3}$. Since $\psi_j$ is complex anti-linear in a
neighborhood of $p_j$, it induces a diffeomorphism $\phi_j:M\to M$.
Assuming that the neighborhoods in question are pairwise
disjoint,  the $\phi_j$'s
commute with each other, and if $S\subset \{ 1, \ldots k\}$
is any subset, we may therefore unambiguously define $\phi_S$
to be the composition of those $\phi_j$'s for which $j\in S$.
Now  $J_S=\phi_{S}^*J$ is an integrable complex
structure on $M$ for each $S\subset \{ 1, \ldots k\}$;
for example, $J_\emptyset =J$.

 Let $c_1(X)$ denote the pull-back to $M$ of the first Chern
class of
 $X$ via the blowing-down map $M\to X$, and let $E_1$,
\ldots , $E_k$ be the Poincar\'e duals of the exceptional
divisors corresponding to $p_1$, \ldots , $p_k$. The
first Chern class of $T^{1,0}_{J_S}M$ is then
$$c_1(M,J_S)= c_1(X) +\sum \epsilon_j ,$$
where
$$\epsilon_j=\left\{ \begin{array}{cc}
E_j& \mbox{if } j\in S\\
-E_j&\mbox{if } j\not\in S.
\end{array}\right.$$
 If $g$ is any Riemannian metric on $M$,
the projection of $c_1(M,J_S)$ into the space
$H^+(g)\subset H^2(M,{\Bbb R})$
of self-dual harmonic 2-forms is therefore
$c_1(M,J_S)^+= c_1(X)^+  +\sum \epsilon_j^+.$
However, $c_1(X)^2=c_1^2(X) > 0$, so $ c_1(X)^+\neq 0$. Now choose
$S$ so that
$$c_1(X)^+\cdot \epsilon_j^+ \geq 0.$$
If $c$ is the  spin$^c$ structure   associated with this choice of
$S$, the
Seiberg-Witten invariant of $(M,H^+,c)$ is non-zero \cite{FM}, and
Theorem
\ref{est}   tells us that
\bea
\frac{1}{32\pi^2}\int_Ms^2d\mu &\geq  &(c_1(M,J_S)^+)^2\\&=&
(c_1(X)^++\sum \epsilon_j^+)^2\\&=&(c_1(X)^+)^2+ 2\sum (c_1(X)^+\cdot
\epsilon_j^+ )+
(\sum \epsilon_j^+)^2\\
&\geq &(c_1(X)^+)^2\\
&\geq & (c_1(X))^2 = c_1^2(X)\\
&=& 2\chi + 3\tau +k
\eea
because the intersection form is positive definite on $H^+=H^+(g)$.

Now suppose  we have a metric $g$ for which
this inequality is actually an equality.
Then each of the inequalities in the above calculation is
 an equality, and
Theorem \ref{est}, applied to the first of these,
tells us that   $g$
is K\"ahler with respect to a  complex
structure $J_g$ compatible with $c$, and hence
satisfying
 $c_1(M,J_g)= c_1(M, J_S)$. By the same reasoning,
$ (\sum\epsilon_j^+)^2=0$, and hence $\sum\epsilon_j^+=0$.
In particular, $c_1(M,J_S)^+=c_1(M,J_{\tilde{S}})^+$, where
$\tilde{S}=\{1, \ldots , k\}-S$, so, even if $b^+=1$, the
Seiberg-Witten invariant of $(M,H^+,\tilde{c})$ is also
non-zero, where $\tilde{c}$ is the spin$^c$ structure
determined by $J_{\tilde{S}}$. The Seiberg-Witten equations
for $\tilde{c}$,
written with respect to the K\"ahler metric $g$,
therefore have an irreducible solution; but this
says \cite{witten,FM} that $-\sum \epsilon_j$
represents
 an effective divisor
on $(M,J_g)$. Thus, if $[\omega ]$
is the K\"ahler class of $(M,g, J_g)$, we have
$[\omega ]\cdot (-\sum \epsilon_j ) > 0$, since
this expression represents  the area of a non-empty holomorphic
curve.
On the other hand, the K\"ahler form $\omega$ is
self-dual with respect to $g$, so we have
$[\omega ]= [\omega ]^+$; thus
 $[\omega ]\cdot \sum \epsilon_j=
 [\omega ]^+\cdot \sum \epsilon_j^+=0$,
in contradiction to the previous assertion.
   Our assumption was therefore
false;  the inequality  is always
  strict.
\end{proof}

\section{Einstein Metrics}

\begin{thm}\label{nein}
Let $X$ be a minimal complex algebraic surface
of general type,
and let $M= X\#k\overline{\Bbb CP}_2$ be obtained from
$X$ by blowing up $k >0$ points.
If $k\geq \frac{2}{3}  c_1^2(X)$,
then $M$ does not admit  Einstein metrics.
 \end{thm}
\begin{proof}
  For any Riemannian
metric $g$ on $M$, one has  the generalized
Gauss-Bonnet formula
$$2\chi + 3\tau =\frac{1}{4\pi^2}\int_M \left(2|W_+|^2 +
\frac{s^2}{24}-\frac{|r_0|^2}{2}\right)d\mu$$
where   $s$,   $r_0$,  and $W_+$  are respectively the
scalar, trace-free Ricci, and  self-dual Weyl
curvatures of $g$;  pointwise norms are
calculated with respect to the metric, and   $d\mu$ is the metric
volume form.
If $g$ is an Einstein metric, $r_0=0$ and
Theorem \ref{best} therefore
implies that
\bea
c_1^2(X)-k = 2\chi+3\tau&=&
 \frac{1}{4\pi^2}\int_M \left(2|W_+|^2+
\frac{s^2}{24}\right)d\mu\\ &>&
\frac{32\pi^2}{4\cdot 24\pi^2}(2 {\chi}+ 3 {\tau}+k)\\
&=&\frac{1}{3}c_1^2(X) ,
\eea
so that
$$\frac{2}{3}  c_1^2(X) > k , $$
contradicting our assumption. Hence $M$ cannot admit
an Einstein metric.
\end{proof}

Our main result  now follows.

\setcounter{main}{0}
\begin{main}
There  are  infinitely many  compact simply-connected
orientable smooth 4-manifolds which
do not
admit Einstein metrics,
but nevertheless
 satisfy the strict
Hitchin-Thorpe inequality.
\end{main}

\begin{proof}
If $X$ is any minimal complex surface of general type with $c_1^2 \geq
3$,
there is then at least one integer
 $k$  satisfying $c_1^2 > k\geq \frac{2}{3}c_1^2$.
The complex surface  $M=X\# k\overline{\Bbb CP}_2$ then
satisfies  the strict Hitchin-Thorpe
inequality
$2\chi > 3|\tau |$,   but does not admit
Einstein metrics by Theorem \ref{nein}.

Now Seiberg-Witten theory implies \cite{FM} that
$c_1^2(X)$ is a diffeomorphism invariant of
 $M=X\# k\overline{\Bbb CP}_2$, so it suffices to
produce a sequence of   simply-connected
minimal surfaces $X_j$ of general type such that
   the sequence of integers $c_1^2(X_j)$ is  increasing.
One such sequence  is given by the
 Fermat surfaces $w^m+x^m+y^m+z^m=0$
of degree $m=j+4$,  with $c_1^2=j^3+4j^2$.
\end{proof}

\section{The Symplectic Case}

In order to keep our discussion as
concrete and elementary as possible, we have
thus far assumed that our 4-manifolds
arose as compact complex surfaces.
The proof of Theorem \ref{best}, however,
 only depends
on the non-vanishing of certain Seiberg-Witten
invariants of $M= X\#k\overline{\Bbb CP}_2$.
Now if  $X$ admits a symplectic structure, the symplectic blow-up
construction of McDuff \cite{dusa}   supplies a
family of such structures on $M$, and
 a result of Taubes \cite{taubes} then provides
us with the non-vanishing invariants we need to prove the
following:

\begin{thm}\label{zest}
Let $(X,\omega )$ be a  symplectic manifold,
and let $M= X\#k\overline{\Bbb CP}_2$.
If $b^+(X)=1$,   assume that   $c_1(X)\cdot [\omega ] < 0$.
Then any Riemannian metric on $M$ satisfies
$$\int_M s^2 d\mu > 32\pi^2  c_1^2(X).$$
 \end{thm}

Here, of course, $c_1(X)$ is the first Chern class of
an almost-complex structure   adapted
to the symplectic structure. The assumption that
 $c_1(X)\cdot [\omega ] < 0$ if
  $b^+=1$ is needed  to
compensate for the fact that
Taubes' proof involves large perturbations
of the Seiberg-Witten equations, whereas the relevant scalar
curvature
estimates stem from the unperturbed equations.

This immediately yields a generalization  of Theorem \ref{nein}:

\begin{thm}\label{nine}
Let $(X,\omega )$ be a  symplectic manifold,
and let $M= X\#k\overline{\Bbb CP}_2$.
If $b^+(X)=1$,   assume that   $c_1(X)\cdot [\omega ] < 0$.
If $k\geq \frac{2}{3}  c_1^2(X)$,
then $M$ does not admit  Einstein metrics.
\end{thm}

Of course, this is a trivial consequence of the Hitchin-Thorpe
inequality unless
  $c_1^2(X) > 0$. On the other hand, it is unnecessarily weak
if $X$ is itself the blow-up of another symplectic manifold.
 In analogy with the Enriques-Kodaira classification, it is
therefore natural to introduce a definition
which characterizes the natural setting for applications
of these results:

\begin{defn}
A minimal symplectic 4-manifold $(X,\omega )$ is  of
{\em    general type}
if
\begin{description}
\item{(a)} $c_1^2(X) > 0$; and
\item{(b)} $c_1(X)\cdot [\omega ] < 0$.
\end{description}
A symplectic 4-manifold   of general type
is   an iterated   symplectic
blow-up of   a  minimal symplectic manifold of general type.
\end{defn}
If $b^+> 1$,
 Taubes \cite{taubes2,taubes3} has shown that   condition (b)
is automatic and that (a) fails only for minimal symplectic
manifolds with $c_1^2=0$.  In analogy to the Kodaira classification,
the latter class of minimal symplectic manifolds might be conjectured
to all
arise as elliptic fibrations unless $c_1$ is a torsion class.
\section{Concluding Remarks}

Theorem \ref{est} tells us that any Riemannian metric on
a non-minimal complex surface $M=X\# k\overline{\Bbb CP}_2$
of general type satisfies the scalar-curvature estimate
$$\int_M s^2 d\mu > 32\pi^2 c_1^2(X), $$
where $X$ is the minimal model for $M$.
In fact, this estimate is sharp, at least if $X$
does not contain any $(-2)$-curves. Indeed,
this last assumption implies \cite{bpv} that $c_1(X) < 0$;
thus   $X$ admits \cite{aubin,yau} a K\"ahler-Einstein
metric $\check{g}$, and one then has
$$\int_X s^2_{\check{g}}d\mu_{\check{g}} = 32\pi^2 c_1^2(X) .$$
Let $p_1$, \ldots , $p_k\in X$ be distinct points which will
be blown up to obtain a smooth  model for $M$, and
choose disjoint  complex coordinate charts
 centered on these points   so that
$$\check{g}= \delta + O(\varrho^2)$$
where $\delta$ and $\varrho$ are respectively the Euclidean
metric and radius associated with the chart.
Define $h_1=\delta - \check{g}$. Let
$h_2$ denote the pull-back of the Fubini-Study
metric on ${\Bbb CP}_1$ to ${\Bbb C}^2-0$
via the tautological projection, and
use these same charts to transplant
$h_2$  to a punctured neighborhood
of each of the $p_1$, \ldots , $p_k$.
Let
$\phi: {\Bbb R}\to {\Bbb R}$ be a non-negative smooth function
which is identically $1$ on $(-\infty , \frac{1}{2})$ and
identically $0$ on $(1,\infty ) $, and, for each sufficiently small $t
> 0$,
let $g_t$ be the smooth Riemannian metric on the blow-up $M$
whose restriction to the open dense set
$X-\{ p_1, \ldots , p_k\}$
is  given by
$$g_t=\check{g} + \phi (t^{-1}\varrho) (h_1+t^4h_2) .$$
For $\varrho < t/2$, this metric coincides up to scale with the
Burns metric \cite{leb0} on the blow-up of ${\Bbb C}^2$ at the origin,
and so has scalar curvature
$s\equiv 0$; and for $\varrho > t$, it coincides with $\check{g}$.
In the transition region $\varrho \in (t/2 , t)$,
one has
\bea \|\phi (t^{-1}\varrho) (t^4h_2-h_1)\| &\leq C_0 t^{2}\\
 \| {\Bbb D}\left[ \phi (t^{-1}\varrho) (t^4h_2-h_1) \right]\| &\leq C_1 t \\
 \| {\Bbb D}^{(2)}\left[ \phi (t^{-1}\varrho) (t^4h_2-h_1) \right]
\| &\leq C_2\eea
where ${\Bbb D}$  is the  Euclidean derivative operator   associated
with the given coordinate system,  and the constants $C_j$ are
independent
of $t$. Thus  $s^2(g_t)$ is uniformly bounded as $t\to 0$,
and since the volume of the annular  transition region
is of order $t^4$, we conclude that
$$\lim_{t\to 0^+} \int_M s^2_{g_t}d\mu_{g_t} =\int_X
s^2_{\check{g}}d\mu_{\check{g}}
= 32\pi^2 c_1^2(X) .$$
The bound is therefore sharp, as claimed.

Even if $X$ contains $(-2)$-curves, the above conclusion
should still hold. Indeed, if
 $\check{X}$ is the complex orbifold obtained by
collapsing all the $(-2)$-curves in $X$, then the Aubin-Yau proof,
without
essential  alterations,  seems to show that
$\check{X}$
admits an  orbifold K\"ahler-Einstein metric $\check{g}$
with   singularities modeled  on ${\Bbb C}^2/(\pm 1)$.
The previous formula for $g_t$ then merely needs to be augmented
by using     Eguchi-Hanson metrics   to smooth out the
orbifold singularities of $\check{g}$.

One  strategy for finding Einstein metrics on a
compact 4-manifold   is to try to minimize
  $\int s^2 d\mu$, since a
critical point of this functional is either
Einstein or scalar-flat. We have just seen, however, that
this strategy will fail on a non-minimal surface of
general type because a minimizing sequence can
collapse (bubble off) to a metric on a topologically
simpler manifold.
In light of this and Theorem \ref{nein}, it seems
plausible to conjecture that non-minimal surfaces
of general type never admit Einstein metrics.
Indeed, making a considerable leap of faith beyond
\cite{leb}, one might  conjecture that Einstein metrics
on irrational surfaces are always K\"ahler.
Further progress in this direction, however, would seem
to require estimates on the norm of the Weyl
curvature which for the present remain  elusive.

\end{document}